\documentclass{article}%
\usepackage{amsmath}
\usepackage{amsfonts}
\usepackage{amssymb}
\usepackage{graphicx}%
\providecommand{\U}[1]{\protect\rule{.1in}{.1in}}
\newtheorem{theorem}{Theorem} [section]

\newtheorem{conjecture}[theorem]{Conjecture}
\newtheorem{corollary}[theorem]{Corollary}

\newtheorem{lemma}[theorem]{Lemma}

\newtheorem{proposition}[theorem]{Proposition}
\newtheorem{remark}[theorem]{Remark}

\newenvironment{proof}[1][Proof]{\noindent\textbf{#1.} }{\ \rule{0.5em}{0.5em}}
\begin{document}

\title{The independence polynomial of a graph at $-1$}
\author{Vadim E. Levit\\Department of Computer Science and Mathematics\\Ariel University Center of Samaria, Ariel, Israel\\levitv@ariel.ac.il
\and Eugen Mandrescu\\Department of Computer Science\\Holon Institute of Technology, Holon, Israel\\eugen\_m@hit.ac.il}
\date{}
\maketitle

\begin{abstract}
The \textit{stability number }$\alpha(G)$ of the graph $G$ is the size of a
maximum stable set of $G$. If $s_{k}$ denotes the number of stable sets of
cardinality $k$ in graph $G$, then $I(G;x)=s_{0}+s_{1}x+...+s_{\alpha
}x^{\alpha}$ is the \textit{independence polynomial} of $G$ \cite{GutHar},
where $\alpha=\alpha(G)$ is the size of a maximum stable set.

In this paper we prove that $I(G;-1)$ satisfies\textit{\ }$\left\vert
I(G;-1)\right\vert \leq2^{\nu(G)}$, where $\nu(G)$ equals the cyclomatic
number of $G$, and the bounds are sharp.

In particular, if $G$ is a connected well-covered graph of girth $\geq6$,
non-isomorphic to $C_{7}$ or $K_{2}$ (e.g., a well-covered tree $\neq K_{2}$),
then $I(G;-1)=0$.

\textbf{Keywords:} stable set, independence polynomial, cyclomatic number,
tree, well-covered graph.

\end{abstract}

\section{Introduction}

Throughout this paper $G=(V,E)$ is a simple (i.e., a finite, undirected,
loopless and without multiple edges) graph with vertex set $V=V(G)$ and edge
set $E=E(G).$ If $X\subset V$, then $G[X]$ is the subgraph of $G$ spanned by
$X$. By $G-W$ we mean the subgraph $G[V-W]$, if $W\subset V(G)$. We also
denote by $G-F$ the partial subgraph of $G$ obtained by deleting the edges of
$F$, for $F\subset E(G)$, and we write shortly $G-e$, whenever $F$ $=\{e\} $.
The \textit{neighborhood} of a vertex $v\in V$ is the set $N_{G}(v)=\{w:w\in
V$ \ \textit{and} $vw\in E\}$, and $N_{G}[v]=N_{G}(v)\cup\{v\}$; if there is
ambiguity on $G$, we use $N(v)$ and $N[v]$, respectively. A vertex $v$ is
\textit{pendant} if its neighborhood contains only one vertex; an edge $e=uv$
is \textit{pendant} if one of its endpoints is a pendant vertex. $K_{n}%
,P_{n},C_{n},K_{n_{1},n_{2},...,n_{p}}$ denote respectively, the complete
graph on $n\geq1$ vertices, the chordless path on $n\geq1$ vertices, the
chordless cycle on $n\geq3$ vertices, and the complete multipartite graph on
$n_{1}+n_{2}+...+n_{p}$ vertices.

The \textit{disjoint union} of the graphs $G_{1},G_{2}$ is the graph
$G=G_{1}\cup G_{2}$ having as a vertex set the disjoint union of
$V(G_{1}),V(G_{2})$, and as an edge set the disjoint union of $E(G_{1}%
),E(G_{2})$. In particular, $nG$ denotes the disjoint union of $n>1$ copies of
the graph $G$. 

If $G_{1},G_{2}$ are disjoint graphs, then their \textit{Zykov
sum} is the graph $G_{1}+G_{2}$ with $V(G_{1})\cup V(G_{2})$ as a vertex set
and $E(G_{1})\cup E(G_{2})\cup\{v_{1}v_{2}:v_{1}\in V(G_{1}),v_{2}\in
V(G_{2})\}$ as an edge set.

As usual, a \textit{tree} is an acyclic connected graph.

A \textit{stable} set in $G$ is a set of pairwise non-adjacent vertices. A
stable set of maximum size will be referred to as a \textit{maximum stable
set} of $G$, and the \textit{stability number }of $G$, denoted by $\alpha(G)
$, is the cardinality of a maximum stable set in $G$.

A graph $G$ is called \textit{well-covered} if all its maximal stable sets are
of the same cardinality \cite{Plum}. If, in addition, $G$ has no isolated
vertices and its order equals $2\alpha(G)$, then $G$ is \textit{very
well-covered} \cite{Fav}. For instance, the graph $G=H\circ K_{1}$, obtained
from $H$ by appending a single pendant edge to each vertex of $H$, is very
well-covered and $\alpha(G)=\left\vert V(H)\right\vert $. The following result
shows that, under certain conditions, any well-covered graph has this form.

\begin{theorem}
\cite{FinHarNow}\label{th3} Let $G$ be a connected graph of girth $\geq6$,
which is isomorphic to neither $C_{7}$ nor $K_{1}$. Then $G$ is well-covered
if and only if $G=H\circ K_{1}$, for some graph $H$ of girth $\geq6$.
\end{theorem}

In other words, Theorem \ref{th3} shows that apart from $K_{1}$ and $C_{7}$,
connected well-covered graphs of girth $\geq6$ are very well-covered.

\begin{proposition}
\cite{Ravindra}\label{prop6} A tree $T$ is well-covered if and only if either
$T=K_{1}$ or $T=H\circ K_{1}$ for some tree $H$.
\end{proposition}

The structure of very well-covered graphs of girth at least $5$ is described
by the following theorem.

\begin{theorem}
\cite{LevMan2007}\label{th1} Let $G$ be a graph of girth at least $5$. Then
$G$ is very well-covered if and only if $G=H\circ K_{1}$, for some graph $H$
of girth $\geq5$.
\end{theorem}

Notice that a well-covered graph can have non-well-covered subgraphs; e.g.,
each subgraph of $C_{5}$ isomorphic to $P_{3}$ is not well-covered, while
$C_{5}$ is well-covered.

\begin{proposition}
\cite{Campbell}\label{prop5} If $G$ is a non-complete well-covered graph, then
$G-N[v]$ is well-covered for any $v\in V(G)$.
\end{proposition}

Let $s_{k}$ be the number of stable sets in $G$ of cardinality $k,k\in
\{1,...,\alpha(G)\}$. The polynomial%
\[
I(G;x)=s_{0}+s_{1}x+s_{2}x^{2}+...+s_{\alpha}x^{\alpha},\ \alpha=\alpha(G),
\]
is called the \textit{independence polynomial} of $G$ \cite{GutHar}. Some
properties of the independence polynomial are presented in
\cite{AlMalSchErdos,Brown,HoedeLi,LeMan04,LeMa04b,LeMa06,LevMan2007,LeMa04e}.
As examples, we mention that:%
\begin{align*}
I(G_{1}\cup G_{2};x) &  =I(G_{1};x)\cdot I(G_{2};x),\\
I(G_{1}+G_{2};x) &  =I(G_{1};x)+I(G_{2};x)-1.
\end{align*}

The following result provides an easy recursive technique in evaluating
independence polynomials of various graphs.

\begin{proposition}
\cite{GutHar,HoedeLi}\label{prop1} If $w\in V(G)$ and $uv\in E(G)$,
then the following equalities hold:

\emph{(i)} $I(G;x)=I(G-w;x)+x\cdot I(G-N[w];x)$;

\emph{(ii)} $I(G;x)=I(G-uv;x)-x^{2}\cdot I(G-N(u)\cup N(v);x)$.
\end{proposition}

The value of a graph polynomial at a specific point can give sometimes a very
surprising information about the structure of the graph (see, for instance,
\cite{BalBolCutPeb}, where the value of the so-called interlace polynomial at
$-1$ is involved). In the case of independence polynomials, let us notice that
if $I(G;x)=s_{0}+s_{1}x+s_{2}x^{2}+...+s_{\alpha}x^{\alpha},\ \alpha
=\alpha(G)$, then:

\begin{itemize}
\item $I(G;1)=s_{0}+s_{1}+s_{2}+...+s_{\alpha}=f(G)$ equals the number of
stable sets of $G$, where $f(G)$ is known as the \textit{Fibonacci number} of
$G$ \cite{KTiWagZieg,PedVest,ProdTichy};

\item $I(G;-1)=s_{0}-s_{1}+s_{2}-...+(-1)^{\alpha}s_{\alpha}=f_{0}%
(G)-f_{1}(G)$, where
\[
f_{0}(G)=s_{0}+s_{2}+s_{4}+...,\quad f_{1}(G)=s_{1}+s_{3}+s_{4}+...
\]
are equal to the numbers of stable sets of even size and odd size of $G$,
respectively. $I(G;-1)$\ is known as the \textit{alternating number of
independent sets} \cite{BSN07}.
\end{itemize}

The difference $\left\vert f_{0}(G)-f_{1}(G)\right\vert $ can be indefinitely
large. It is easy to check that the complete $n$-partite graph $K_{\alpha
,\alpha,...,\alpha}$ is well-covered, $\alpha(K_{\alpha,\alpha,...,\alpha
})=\alpha$, and its independence polynomial is $I(K_{\alpha,\alpha,...,\alpha
};x)=n(1+x)^{\alpha}-(\alpha-1)$. Hence, $I(K_{\alpha,\alpha,...,\alpha
};-1)=1-\alpha$, i.e., for any negative integer $k$ there is some connected
well-covered graph $G$ such that $I(G;-1)=k$.

Let $G_{i},1\leq i\leq k$, be $k$ graphs with $I(G_{i};-1)=2$ for every
$i\in\{1,2,...,k\}$, and $H=G_{1}+G_{2}+...+G_{k}$. Then,
\[
I(H;x)=I(G_{1};x)+I(G_{2};x)+...+I(G_{k};x)-(k-1)
\]
and consequently, $I(H;-1)=2k-(k-1)=k+1$. In other words, for any positive
integer $k$ there is some connected well-covered graph $G$ such that
$I(G;-1)=k$.

In this paper we prove that:

\begin{itemize}
\item $I(T;-1)\in\{-1,0,1\}$ for every tree $T$;

\item $I(G;-1)=0$ for every connected well-covered graph $G$ of girth $\geq6$,

non-isomorphic to $C_{7}$ or $K_{2}$;

\item $\left\vert I(G;-1)\right\vert \leq2^{\nu(G)}$, for every graph $G$,
where $\nu(G)$ is its cyclomatic number.
\end{itemize}

\section{$I(G;-1)$ for a graph $G$ with at most one cycle}

There are graphs $G$ with at least one pendant vertex having $I(G;-1)\in
\{-1,0,1\}$; see, for instance, the graphs from Figure \ref{fig44}, whose
independence polynomials are, respectively,%
\begin{align*}
I(G_{1};x)  &  =1+5x+5x^{2}+x^{3},I(G_{2};x)=\left(  1+5x+4x^{2}+x^{3}\right)
\left(  1+x\right)  ,\\
I(G_{3};x)  &  =1+4x+2x^{2},I(G_{4};x)=1+5x+5x^{2}+2x^{3}.
\end{align*}
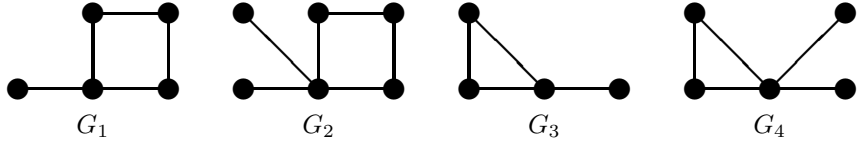
\begin{figure}[h]
\setlength{\unitlength}{1cm}\begin{picture}(5,1.4)\thicklines
\multiput(1,0.5)(1,0){3}{\circle*{0.29}}
\multiput(2,1.5)(1,0){2}{\circle*{0.29}}
\put(1,0.5){\line(1,0){2}}
\put(2,0.5){\line(0,1){1}}
\put(2,1.5){\line(1,0){1}}
\put(3,0.5){\line(0,1){1}}
\put(2,0){\makebox(0,0){$G_{1}$}}
\multiput(4,0.5)(1,0){3}{\circle*{0.29}}
\multiput(4,1.5)(1,0){3}{\circle*{0.29}}
\put(4,0.5){\line(1,0){2}}
\put(5,1.5){\line(1,0){1}}
\put(6,0.5){\line(0,1){1}}
\put(5,0.5){\line(0,1){1}}
\put(5,0.5){\line(-1,1){1}}
\put(5,0){\makebox(0,0){$G_{2}$}}
\multiput(7,0.5)(1,0){3}{\circle*{0.29}}
\put(7,1.5){\circle*{0.29}}
\put(7,0.5){\line(1,0){2}}
\put(7,0.5){\line(0,1){1}}
\put(7,1.5){\line(1,-1){1}}
\put(8,0){\makebox(0,0){$G_{3}$}}
\multiput(10,0.5)(1,0){3}{\circle*{0.29}}
\multiput(10,1.5)(2,0){2}{\circle*{0.29}}
\put(10,0.5){\line(1,0){2}}
\put(10,0.5){\line(0,1){1}}
\put(10,1.5){\line(1,-1){1}}
\put(11,0.5){\line(1,1){1}}
\put(11,0){\makebox(0,0){$G_{4}$}}
\end{picture}
\caption{$I(G_{1};-1)=0=I(G_{2};-1)$, while $I(G_{3};-1)=-1=-I(G_{4};-1)$.}%
\label{fig44}%
\end{figure}

\begin{lemma}
\label{lem2}If $u\in V(G)$ is a pendant vertex of $G$ and $v\in N(u)$, then
\[
I(G;-1)=(-1)\cdot I(G-N[v];-1).
\]
Moreover, if $G$ has two pendant vertices at $3$ distance apart, then
$I(G;-1)=0$.
\end{lemma}

\begin{proof}
Since $u\in V(G)$ is a pendant vertex of $G$ and $v\in N(u)$, Proposition
\ref{prop1}\emph{(i)} assures that
\begin{align*}
I(G;x)  &  =I(G-v;x)+x\cdot I(G-N[v];x)\\
&  =(1+x)\cdot I(G-\{u,v\};x)+x\cdot I(G-N[v];x)
\end{align*}
and this implies $I(G;-1)=(-1)\cdot I(G-N[v];-1)$.

Let $a,b$ be two pendant vertices of $G$ with $dist(a,b)=3$ and let $v\in
N(a)$. According to Proposition \ref{prop1}\emph{(i)}, we get:%
\begin{align*}
I(G;x)  &  =I(G-v;x)+x\cdot I(G-N[v];x)\\
&  =I(\{a\};x)\cdot I(G-\{a,v\};x)+x\cdot I(\{b\};x)I(G-\{b\}\cup N[v];x)\\
&  =(1+x)\left[  I(G-\{a,v\};x)+x\cdot I(G-\{b\}\cup N[v];x)\right]  ,
\end{align*}
which clearly implies $I(G;-1)=0$.
\end{proof}

\begin{corollary}
\label{cor3}If $G$ is a well-covered graph of girth $\geq6$ and $G\neq qC_{7}$
for any $q\geq1$, or $G$ is a very well-covered graph of girth at least $5$
and $G\neq qK_{2}$ for any $q\geq1$, then $I(G;-1)=0$, i.e., the number of
stable sets of even size equals the number of stable sets of odd size. In
particular, the assertion is true for every well-covered tree $T\neq K_{2}$.
\end{corollary}

\begin{proof}
Notice that $I(K_{2};x)=1+2x$, and $I(C_{7};x)=1+7x+14x^{2}+7x^{3}$.
Therefore, $I(K_{2};-1)=-1$, while $I(C_{7};-1)=1$.

If $G=K_{1}$, then $I(K_{1};x)=1+x$ and, clearly, $I(G;-1)=0$.

Otherwise, according to Theorems \ref{th3} and \ref{th1}, it follows that
$G=H\circ K_{1}$, for some graph $H$ of order at least two, since $G\neq K_{2}
$. Consequently, $G$ has at least two pendant vertices at $3$ distance apart,
and by Lemma \ref{lem2}, we obtain that $I(G;-1)=0$.

If $T\neq K_{2}$ is a well-covered tree, then either $T=K_{1}$ or $T$ is very
well-covered and its girth is greater than $5$. In both cases, we get
$I(T;-1)=0$.
\end{proof}

In \cite{Arocha} it is shown that
\[
I(P_{n};x)=F_{n+1}(x)\text{ and }I(C_{n},x)=F_{n-1}(x)+2x\cdot F_{n-2}(x),
\]
where $F_{n}(x),n\geq0,$ are \textit{Fibonacci polynomials}, i.e., the
polynomials defined recursively by%
\[
F_{0}(x)=1,F_{1}(x)=1,F_{n}(x)=F_{n-1}(x)+x\cdot F_{n-2}(x).
\]
Let us notice that $P_{n},n\geq5$, and $C_{n}$, for $n=6$ or $n\geq8$, are not well-covered.

\begin{lemma}
\label{lem1}For $n\geq1$, the following equalities hold:

\emph{(i)} $I(P_{3n-2};-1)=0$ and $I(P_{3n-1};-1)=I(P_{3n};-1)=(-1)^{n}$;

\emph{(ii)} $I(C_{3n};-1)=2\cdot(-1)^{n},I(C_{3n+1};-1)=(-1)^{n}$ and
$I(C_{3n+2};-1)=(-1)^{n+1}$.
\end{lemma}

\begin{proof}
\emph{(i)} We prove by induction on $n$.

For $n=1$, we have $I(P_{1};-1)=0$ and $I(P_{2};-1)=I(P_{3};-1)=-1$, because
\[
I(P_{1};x)=1+x,I(P_{2};x)=1+2x,I(P_{3};x)=1+3x+x^{2}.
\]

Assume that the assertion is true for any $k\leq3n$. Using Proposition
\ref{prop1}\emph{(i)}, we obtain $I(P_{k+1};x)=I(P_{k};x)+xI(P_{k-1};x)$,
which leads respectively to:%
\begin{align*}
I(P_{3n+1};-1)  &  =I(P_{3n};-1)-I(P_{3n-1};-1)=(-1)^{n}-(-1)^{n}=0;\\
I(P_{3n+2};-1)  &  =I(P_{3n+1};-1)-I(P_{3n};-1)=0-(-1)^{n}=(-1)^{n+1};\\
I(P_{3n+3};-1)  &  =I(P_{3n+2};-1)-I(P_{3n+1};-1)=(-1)^{n+1}-0=(-1)^{n+1}.
\end{align*}

\emph{(ii)} Firstly,
\begin{align*}
I(C_{3};x)  &  =1+3x,\quad I(C_{4};x)=1+4x+2x^{2},\\
I(C_{5};x)  &  =1+5x+5x^{2},\quad I(C_{6};x)=1+6x+9x^{2}+2x^{3},
\end{align*}
which implies that
\begin{align*}
I(C_{3};-1)  &  =2\cdot(-1),\quad I(C_{3+1};-1)=(-1),\\
I(C_{3+2};-1)  &  =(-1)^{2},\quad I(C_{3\cdot2};-1)=2\cdot(-1)^{2}.
\end{align*}
Let $uv$ be an edge of some $C_{n},n\geq7$. By Proposition \ref{prop1}%
\emph{(i)}, we deduce that%
\begin{align*}
I(P_{n};x)  &  =I(P_{n-1};x)+x\cdot I(P_{n-2};x)=...\\
&  =(1+2x)\cdot I(P_{n-3};x)+x\cdot(1+x)\cdot I(P_{n-4};x).
\end{align*}
Now, using Proposition \ref{prop1}\emph{(ii)}, we get%
\begin{align*}
I(C_{n};x)  &  =I(C_{n}-uv;x)-x^{2}\cdot I(C_{n}-N(u)\cup N(v);x)\\
&  =I(P_{n};x)-x^{2}\cdot I(P_{n-4};x)=(1+2x)\cdot I(P_{n-3};x)+x\cdot
I(P_{n-4};x).
\end{align*}
Hence, we obtain $I(C_{n};-1)=(-1)\cdot\lbrack I(P_{n-3};-1)+I(P_{n-4};-1)]$.
Since, by part \emph{(i)},
\[
I(P_{n-3};-1)+I(P_{n-4};-1)\in\{0+(-1)^{k},2\cdot(-1)^{k}\},
\]
where $k$ depends on $n-3\in\{3k-2,3k-1,3k\}$, it is easy to get that
\[
I(C_{3k};-1)=2\cdot(-1)^{k},\quad I(C_{3k+1};-1)=(-1)^{k}\text{ and
}I(C_{3k+2};-1)=(-1)^{k+1},
\]
and this completes the proof.
\end{proof}

\begin{figure}[h]
\setlength{\unitlength}{1cm}\begin{picture}(5,1.2)\thicklines
\multiput(1.5,0)(1,0){3}{\circle*{0.29}}
\put(2.5,1){\circle*{0.29}}
\put(1.5,0){\line(1,0){2}}
\put(2.5,0){\line(0,1){1}}
\put(1,0.5){\makebox(0,0){$T_{1}$}}
\multiput(5,0)(1,0){4}{\circle*{0.29}}
\put(6,1){\circle*{0.29}}
\put(7,1){\circle*{0.29}}
\put(6,0){\line(0,1){1}}
\put(5,0){\line(1,0){3}}
\put(7,0){\line(0,1){1}}
\put(4.5,0.5){\makebox(0,0){$T_{2}$}}
\multiput(10,0)(1,0){3}{\circle*{0.29}}
\multiput(10,1)(1,0){4}{\circle*{0.29}}
\put(10,0){\line(1,0){2}}
\put(10,0){\line(1,1){1}}
\put(10,0){\line(0,1){1}}
\put(12,0){\line(1,1){1}}
\put(12,0){\line(0,1){1}}
\put(9.5,0.5){\makebox(0,0){$T_{3}$}}
\end{picture}
\caption{Non-well-covered tees.}%
\label{fig707}%
\end{figure}
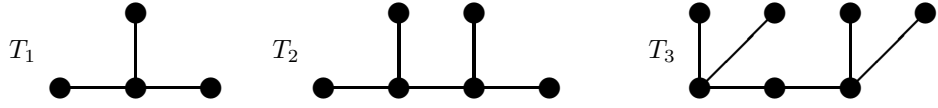

Let us notice that there exist non-well-covered trees $T_{1},T_{2},T_{3}$,
non-isomorphic to $P_{n}$ (see Figure \ref{fig707}), such that $I(T_{1}%
;-1)=-1,I(T_{2};-1)=0,I(T_{3};-1)=1$, because
\begin{align*}
I(T_{1};x)  &  =1+4x+3x^{2}+x^{3},I(T_{2};x)=1+6x+10x^{2}+6x^{3}+x^{4},\\
I(T_{3};x)  &  =1+7x+15x^{2}+12x^{3}+5x^{4}+x^{5}.
\end{align*}

The following theorem generalizes Lemma \ref{lem1}\emph{(i)}.

\begin{theorem}
\label{th6}For any tree $T,I(T;-1)\in\{-1,0,1\}$, i.e., the number of stable
sets of even size varies by at most one from the number of stable sets of odd size.
\end{theorem}

\begin{proof}
Let us observe that it is sufficient to show that $I(T;-1)\in\{-1,0,1\}$,
because $I(T;x)=s_{0}+s_{1}x+...+s_{\alpha}x^{\alpha},\alpha=\alpha(T)$,
implies
\[
I(T;-1)=(s_{0}+s_{2}+s_{4}+...)-(s_{1}+s_{3}+s_{4}+...).
\]
We prove by induction on $n=\left\vert V(T)\right\vert $.

For $n=1,I(T;x)=1+x$ and $I(T;-1)=0$.

For $n=2,I(T;x)=1+2x$ and $I(T;-1)=-1$.

For $n=3,I(T;x)=1+3x+x^{2}$ and $I(T;-1)=-1$.

For $n=4$, only two trees are non-isomorphic, namely $P_{4}$ and $K_{1,3}$.
Since
\[
I(P_{4};x)=1+4x+3x^{2}\text{ and }I(K_{1,3};x)=1+4x+3x^{2}+x^{3},
\]
it follows that $I(P_{4};-1)=0$, while $I(K_{1,3};-1)=-1$.

Finally, for $n=5$ there are three non-isomorphic trees, namely $K_{1,4}%
,T_{5}$ and $P_{5}$ (see Figure \ref{fig717}). We have successively%
\begin{align*}
I(K_{1,4};x)  &  =1+5x+6x^{2}+4x^{3}+x^{4},\\
I(T_{5};x)  &  =1+5x+6x^{2}+2x^{3},\quad I(P_{5};x)=1+5x+6x^{2}+x^{3},
\end{align*}
which give $I(K_{1,4};-1)=-1,I(T_{5};-1)=0$ and $I(P_{5};-1)=1$.
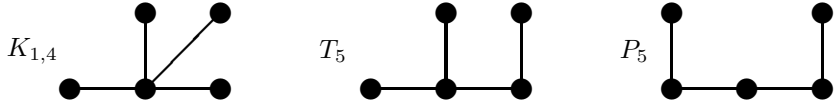
\begin{figure}[h]
\setlength{\unitlength}{1cm}\begin{picture}(5,1.2)\thicklines
\multiput(2,0)(1,0){3}{\circle*{0.29}}
\put(3,1){\circle*{0.29}}
\put(4,1){\circle*{0.29}}
\put(2,0){\line(1,0){2}}
\put(3,0){\line(0,1){1}}
\put(3,0){\line(1,1){1}}
\put(1.5,0.5){\makebox(0,0){$K_{1,4}$}}
\multiput(6,0)(1,0){3}{\circle*{0.29}}
\put(7,1){\circle*{0.29}}
\put(8,1){\circle*{0.29}}
\put(7,0){\line(0,1){1}}
\put(6,0){\line(1,0){2}}
\put(8,0){\line(0,1){1}}
\put(5.5,0.5){\makebox(0,0){$T_{5}$}}
\multiput(10,0)(1,0){3}{\circle*{0.29}}
\multiput(10,1)(2,0){2}{\circle*{0.29}}
\put(10,0){\line(1,0){2}}
\put(10,0){\line(0,1){1}}
\put(12,0){\line(0,1){1}}
\put(9.5,0.5){\makebox(0,0){$P_{5}$}}
\end{picture}
\caption{The three non-isonorphic trees on $5$ vertices.}%
\label{fig717}%
\end{figure}

Let us suppose that $T=(V,E)$ is a tree with $\left\vert V\right\vert
=n+1\geq5,v$ is a pendant vertex of $T,N(v)=\{u\}$ and $T_{1},T_{2},...,T_{m}$
are the trees of the forest $T-N[u]$. According to Proposition \ref{prop1}%
\emph{(ii)}, we get the following:%
\begin{align*}
I(T;x)  &  =I(T-uv;x)-x^{2}\cdot I(T-N(u)\cup N(v);x)=\\
&  =(1+x)\cdot I(T-v;x)-x^{2}\cdot I(T-N[u];x)=\\
&  =(1+x)\cdot I(T-v;x)-x^{2}\cdot I(T_{1};x)\cdot I(T_{2};x)\cdot...\cdot
I(T_{m};x).
\end{align*}
Consequently, $I(T;-1)=-I(T_{1};-1)\cdot I(T_{2};-1)\cdot...\cdot
I(T_{m};-1)\in\{-1,0,1\}$, since every tree $T_{i}$ has less than $n$
vertices, and by the induction hypothesis, $I(T_{i};-1)\in\{-1,0,1\}$.
\end{proof}

\begin{corollary}
\emph{(i)} For every tree, the number of dependent sets of even size varies by
at most one from the number of dependent sets of odd size.

\emph{(ii)} If $T$ is well-covered tree and $T\neq K_{2}$, then the number of
dependent sets of even size equals the number of dependent sets of odd size.

\emph{(iii)} If $T$ is well-covered tree and $T\neq K_{2}$, then the number of
all stable sets and the number of all dependent sets are even.
\end{corollary}

\begin{proof}
Let $p_{1},p_{2},$ be the numbers of stable sets of the tree $T$, of odd and
even size, respectively, and $q_{1},q_{2}$ be the numbers of dependent sets in
$T$, of odd and even size, respectively. Clearly, if $\left\vert
V(T)\right\vert =n$, then $p_{1}+q_{1}=p_{2}+q_{2}=2^{n-1}$. According to
Theorem \ref{th6}, $\left\vert p_{1}-p_{2}\right\vert \leq1$ and this implies
$\left\vert q_{1}-q_{2}\right\vert \leq1$, while for a well-covered tree
$T\neq K_{2},p_{1}=p_{2}$ and this leads to $q_{1}=q_{2}$. Hence for
well-covered trees different from $K_{2}$ both $p_{1}+p_{2}$ and $q_{1}+q_{2}
$ are even numbers.
\end{proof}

\begin{corollary}
\label{cor1}If $F$ is a forest, then $I(F;-1)\in\{-1,0,1\}$, and $I(F;-1)=0$
whenever at least one of its components is a well-covered tree different from
$K_{2}$.
\end{corollary}

\begin{proof}
If $T_{1},T_{2},...,T_{m},m\geq1$, are the connected components of $F$, then
\[
I(F;x)=I(T_{1};x)\cdot I(T_{1};x)\cdot...\cdot I(T_{m};x),
\]
because $F=T_{1}\cup T_{2}\cup...\cup T_{m}$. Now, the conclusions follow from
Theorem \ref{th6} and Corollary \ref{cor3}.
\end{proof}

The unicycle non-well-covered graph $G_{1}$ from Figure \ref{fig23} has
$I(G_{1};x)=1+5x+5x^{2}+x^{3}$ and $I(G_{1};-1)=0$. 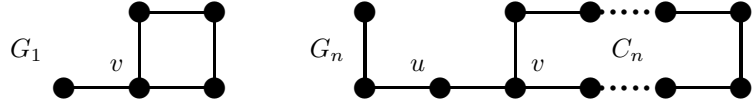
\begin{figure}[h]
\setlength{\unitlength}{1cm}\begin{picture}(5,1)\thicklines
\multiput(2,0)(1,0){3}{\circle*{0.29}}
\multiput(3,1)(1,0){2}{\circle*{0.29}}
\put(2,0){\line(1,0){2}}
\put(3,0){\line(0,1){1}}
\put(3,1){\line(1,0){1}}
\put(4,0){\line(0,1){1}}
\put(2.7,0.3){\makebox(0,0){$v$}}
\put(1.5,0.5){\makebox(0,0){$G_{1}$}}
\multiput(6,0)(1,0){6}{\circle*{0.29}}
\put(6,1){\circle*{0.29}}
\multiput(8,1)(1,0){4}{\circle*{0.29}}
\put(6,0){\line(1,0){3}}
\put(6,0){\line(0,1){1}}
\put(8,0){\line(0,1){1}}
\put(8,1){\line(1,0){1}}
\put(10,0){\line(1,0){1}}
\put(10,1){\line(1,0){1}}
\put(11,0){\line(0,1){1}}
\multiput(9,0)(0.15,0){7}{\circle*{0.07}}
\multiput(9,1)(0.15,0){7}{\circle*{0.07}}
\put(6.7,0.3){\makebox(0,0){$u$}}
\put(8.3,0.3){\makebox(0,0){$v$}}
\put(9.5,0.5){\makebox(0,0){$C_{n}$}}
\put(5.5,0.5){\makebox(0,0){$G_{n}$}}
\end{picture}
\caption{$I(G_{1};-1)=0$, while $I(G_{n};-1)\in\{-2,-1,1,2\}$, for $n\geq4$.}%
\label{fig23}%
\end{figure}

Let us notice that for $n\geq4$, the unicycle graph $G_{n}$ in Figure
\ref{fig23} has
\begin{align*}
I(G_{n};x)  &  =I(G_{n}-uv;x)-x^{2}\cdot I(G_{n}-N(u)\cup N(v);x)\\
&  =I(P_{3};x)\cdot I(C_{n};x)-x^{2}\cdot(1+x)\cdot I(P_{n-3};x),
\end{align*}
and hence, by Lemma \ref{lem1}, we obtain
\[
I(G_{n};-1)=I(P_{3};-1)\cdot I(C_{n};-1)=(-1)\cdot I(C_{n};-1)\in
\{-2,-1,1,2\}.
\]

The case of unicycle well-covered graphs is more specific.

\begin{proposition}
If $G$ is a unicycle well-covered graph and $G\neq C_{3}$, then $I(G;-1)\in
\{-1,0,1\}$.
\end{proposition}

\begin{proof}
Let us notice that $I(C_{3};-1)=-2$, because $I(C_{3};x)=1+3x$.

If $G$ is disconnected, then each cycle-free component $H$ is a well-covered
graph, which, by Corollary \ref{cor3}, contributes with $I(H;-1)\in\{-1,0\}$
in the product that equals $I(G;-1)$. Therefore, we may assume that $G$ is connected.

We show that the assertion is true by induction on $n=\left\vert
V(G)\right\vert $.

If $G\in\{C_{4},C_{5},C_{7}\}$, then Lemma \ref{lem1} ensures that
$I(G;-1)\in\{-1,0,1\}$.

If $G$ contains $C_{3}$ and $G\neq C_{3}$, then for $n=4$\ there is no such
connected graph, while for $n\in\{5,6\}$, the only well-covered graphs are
depicted in Figure \ref{fig112}. In each case, $I(G;-1)\in\{-1,0,1\}$.
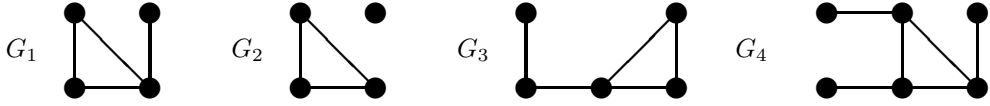
\begin{figure}[h]
\setlength{\unitlength}{1cm}\begin{picture}(5,1)\thicklines
\multiput(1,0)(1,0){2}{\circle*{0.29}}
\multiput(1,1)(1,0){2}{\circle*{0.29}}
\put(1,0){\line(1,0){1}}
\put(1,1){\line(1,-1){1}}
\put(1,0){\line(0,1){1}}
\put(2,0){\line(0,1){1}}
\put(0.3,0.5){\makebox(0,0){$G_{1}$}}
\multiput(4,0)(1,0){2}{\circle*{0.29}}
\multiput(4,1)(1,0){2}{\circle*{0.29}}
\put(4,0){\line(1,0){1}}
\put(4,0){\line(0,1){1}}
\put(4,1){\line(1,-1){1}}
\put(3.3,0.5){\makebox(0,0){$G_{2}$}}
\multiput(7,0)(1,0){3}{\circle*{0.29}}
\multiput(7,1)(2,0){2}{\circle*{0.29}}
\put(7,0){\line(1,0){2}}
\put(7,0){\line(0,1){1}}
\put(8,0){\line(1,1){1}}
\put(9,0){\line(0,1){1}}
\put(6.3,0.5){\makebox(0,0){$G_{3}$}}
\multiput(11,0)(1,0){3}{\circle*{0.29}}
\multiput(11,1)(1,0){3}{\circle*{0.29}}
\put(11,0){\line(1,0){2}}
\put(11,1){\line(1,0){1}}
\put(12,0){\line(0,1){1}}
\put(12,1){\line(1,-1){1}}
\put(13,0){\line(0,1){1}}
\put(10,0.5){\makebox(0,0){$G_{4}$}}
\end{picture}
\caption{All $G_{i},1\leq i\leq4$, have $C_{3}$ as unique cycle, but only
$G_{2},G_{3}$ and $G_{4}$ are well-covered.}%
\label{fig112}%
\end{figure}

Assume that the statement is true for unicycle well-covered graphs having at
most $n$ vertices, and let $G$ be such a graph with $\left\vert
V(G)\right\vert =n+1$. Since $G$ is connected and $G\notin\{C_{3},C_{4}%
,C_{5},C_{7}\}$, we infer that $G$ has at least one pendant vertex. Hence
Lemma \ref{lem2} implies that $I(G;-1)=(-1)\cdot I(G-N[v];-1)$. On the other
hand, Proposition \ref{prop5} assures that $G-N[v]$ is a well-covered graph.
Let $H_{i},1\leq i\leq q$, be all the connected components of $G-N[v]$.

If $G-N[v]$ is forest, then $I(G;-1)\in\{-1,0,1\}$, according to Theorem
\ref{th6}.

If $G-N[v]$ is still a unicycle graph, let $H_{1}$ be the component containing
the cycle. By Theorem \ref{th6}, $I(H_{2};-1)\cdot...\cdot I(H_{q}%
;-1)\in\{-1,0,1\}$, while by induction hypothesis, $I(H_{1};-1)\in\{-1,0,1\}$,
as well. Therefore, we obtain
\begin{align*}
I(G;-1)  &  =(-1)\cdot I(G-N[v];-1)\\
&  =I(H_{1};-1)\cdot I(H_{2};-1)\cdot...\cdot I(H_{q};-1)\in\{-1,0,1\},
\end{align*}
and this completes the proof.
\end{proof}

\begin{corollary}
\label{cor2}If $T=(V,E)$ is well-covered tree of order $\geq6$, then the
following assertions are true:

\emph{(i)} $I(T-v;-1)=I(T-N[v];-1)=0$ for any $v\in V$ such that $T-N[v]\neq
qK_{2},q\geq1$;

\emph{(ii)} $I(T-uv;-1)=0$ holds for any $uv\in E$;

\emph{(iii)} $I(T-N(u)\cup N(v);-1)=0$ holds for any $uv\in E$.
\end{corollary}

\begin{proof}
\emph{(i)} According to Proposition \ref{prop1}\emph{(i)} and Corollary
\ref{cor3}, we obtain
\[
I(T;-1)=I(T-v;-1)+(-1)\cdot I(T-N[v];-1)=0.
\]
Since $T-N[v]\neq qK_{2},q\geq1$, and according to Proposition \ref{prop5},
$T-N[v]$ is a forest consisting of well-covered trees, Corollary \ref{cor1}
ensures that $I(T-N[v];-1)=0$, which further implies that $I(T-v;-1)=0$.

\emph{(ii)} Let $uv\in E(T)$. If $u$ is a pendant vertex, then
\[
I(T-uv;x)=I(\{u\};x)\cdot I(T-u;x)=(1+x)\cdot I(T-u;x),
\]
which assures that $I(T-uv;-1)=0$. If none of $u,v$ is pendant, then $T-uv$
is, by Proposition \ref{prop6}, a disjoint union of two well-covered trees
$T_{1},T_{2}$, and consequently, Corollary \ref{cor3} implies
\[
I(T-uv;-1)=I(T_{1};-1)\cdot I(T_{2};-1)=0,
\]
because at least one of $T_{1},T_{2}$ is non-isomorphic to $K_{2}$.

\emph{(iii)} By Proposition \ref{prop1}\emph{(iii)}, we have
\[
I(T;x)=I(T-uv;x)-x^{2}\cdot I(T-N(u)\cup N(v);x),
\]
which yields $I(T-N(u)\cup N(v);-1)=I(T-uv;-1)-I(T;-1)=0$.
\end{proof}

\begin{figure}[h]
\setlength{\unitlength}{1cm}\begin{picture}(5,1.2)\thicklines
\multiput(4,0)(1,0){6}{\circle*{0.29}}
\multiput(4,1)(1,0){2}{\circle*{0.29}}
\multiput(8,1)(1,0){2}{\circle*{0.29}}
\put(4,0){\line(0,1){1}}
\put(5,0){\line(0,1){1}}
\put(8,0){\line(0,1){1}}
\put(4,0){\line(1,0){5}}
\put(9,0){\line(0,1){1}}
\end{picture}
\caption{The tree $T$ satisfies $I(T;-1)=I(T-v;-1)=0$ for any vertex $v$.}%
\label{fig101}%
\end{figure}
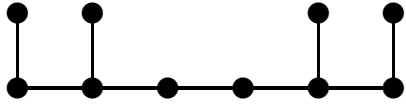

\begin{remark}
There is a non-well-covered tree $T$ satisfying the equalities \emph{(i) -
(iii)} in Corollary \ref{cor2} (see Figure \ref{fig101}).
\end{remark}

\begin{corollary}
For any tree $T$ of order $\geq2$, and for any vertex $v\in V(T)$, it is true
that
\[
I(T-v;-1)\cdot I(T-N[v];-1)\in\{0,1\}\text{ and }I(T;-1)\cdot I(T-v;-1)\in
\{0,1\}.
\]

\end{corollary}

\begin{proof}
According to Proposition \ref{prop1}\emph{(i)}, we get $I(T;x)=I(T-v;x)+x\cdot
I(T-N[v];x)$, which implies that:%
\[
I(T;-1)=I(T-v;-1)-I(T-N[v];-1).
\]
Since $T-v,T-N[v]$ are forests, Corollary \ref{cor1} assures that
\[
I(T;-1),I(T-v;-1),I(T-N[v];-1)\in\{-1,0,1\}
\]
and these are possible provided $I(T-v;-1),I(T-N[v];-1)\in\{0,1\}$.

Further, suppose $w\in V,v\notin V,F=T+v=(V\cup\{v\},E\cup\{vw\})$. According
to Proposition \ref{prop1}\emph{(i)}, we get
\[
I(F;x)=I(F-v;x)+x\cdot I(F-N[v];x)=I(T;x)+x\cdot I(T-w;x),
\]
which implies that:%
\[
I(F;-1)=I(T;-1)-I(T-w;-1).
\]
Since $I(F;-1)\in\{-1,0,1\}$, we infer that if $I(T;-1)=-1$, then
$I(T-v;-1)\neq-1$, and if $I(T;-1)=1$, then $I(T-v;-1)\neq1$. Hence,
$I(T;-1)\cdot I(T-v;-1)\in\{0,1\}$.
\end{proof}

\section{$I(G;-1)$ and the cyclomatic number of $G$}

The \textit{cyclomatic number} $\nu(G)$ of the graph $G$ is the dimension of
the \textit{cycle space} of $G$, i.e., the dimension of the linear space
spanned by the edge sets of all the cycles of $G$. It is known that
$\nu(G)=\left\vert E(G)\right\vert -\left\vert V(G)\right\vert +p$, where $p$
is the number of connected components of $G$. If $e\in E(G)$ belongs to a
cycle, then $G-e$ has the same number of connected components, and hence,
$\nu(G-e)=\left\vert E(G)-\{e\}\right\vert -\left\vert V(G)\right\vert
+p=\nu(G)-1$.

\begin{theorem}
For any graph $G$ the alternating number of independent sets is bounded as
follows
\[
-2^{\nu(G)}\leq I(G;-1)\leq2^{\nu(G)},
\]
where $\nu(G)$ is the cyclomatic number of $G$.
\end{theorem}

\begin{proof}
We prove by induction on $\nu(G)$.

If $\nu(G)=0$, then $G$ is a forest, and according to Corollary \ref{cor1}, we
obtain
\[
I(G;-1)\in\{-1,0,1\}=\{-2^{\nu(G)},0,2^{\nu(G)}\}.
\]

Assume that the assertion is true for graphs with cyclomatic number $\leq k$,
and let $G$ be a graph with $\nu(G)=k+1$. Since $\nu(G)\geq1$, $G$ has at
least one cycle, and if $uv\in E(G)$ belongs to some cycle of $G$, then
$\nu(G-uv)=\nu(G)-1=k$. According to Proposition \ref{prop1}\emph{(ii)}, we
get:
\[
I(G;-1)=I(G-uv;-1)-(-1)^{2}\cdot I(G-N(u)\cup N(v);-1),
\]
which assures that
\[
\left\vert I(G;-1)\right\vert \leq\left\vert I(G-uv;-1)\right\vert +\left\vert
I(G-N(u)\cup N(v);-1)\right\vert \leq2\cdot2^{k},
\]
because $\nu(G-uv)=k$ and $\nu(G-N(u)\cup N(v))\leq k$.
\end{proof}

\begin{remark}
Let $G=qK_{3}$. Then, $I(G;x)=(1+3x)^{q}$ and hence, $I(G;-1)=(-2)^{q}%
=(-2)^{\nu(G)}$.
\end{remark}

\begin{lemma}
\label{lem5}Let $\{G_{i}:1\leq i\leq n\},n>1$, be a family of graphs, $v$ a
vertex belonging to none of $G_{i}$, and $H=H[v,G_{1},...,G_{n}]$ be the graph
obtained by joining $v$ to some vertex $u_{i}\in V(G_{i})$ for every
$i\in\{1,2,...,n\}$. Then
\[
I(H;-1)=I(G_{1};-1)\cdot...\cdot I(G_{n};-1)-I(G_{1}-u_{1};-1)\cdot...\cdot
I(G_{n}-u_{n};-1).
\]
Moreover, $\nu(H)=\nu(G_{1})+...+\nu(G_{n})$ and $H$ is connected whenever
each $G_{i}$ is connected.
\end{lemma}

\begin{proof}
According to Proposition \ref{prop1}\emph{(i)}, we have
\begin{align*}
I(H;x)  &  =I(H-v;x)+x\cdot I(W-N[v];x)\\
&  =I(G_{1};x)\cdot...\cdot I(G_{n};x)+x\cdot I(G_{1}-u_{1};x)\cdot...\cdot
I(G_{n}-u_{n};x),
\end{align*}
which clearly implies
\[
I(H;-1)=I(G_{1};-1)\cdot...\cdot I(G_{n};-1)-I(G_{1}-u_{1};-1)\cdot...\cdot
I(G_{n}-u_{n};-1).
\]

Since each cycle of $H$ appears as a cycle in one of $G_{i}$, it follows that
\[
\nu(H)=\nu(G_{1})+\nu(G_{2})+...+\nu(G_{n}).
\]
Evidently, $H$ is connected whenever every $G_{i}$ is connected.
\end{proof}

\begin{corollary}
\label{cor33}For any graph $G$, there exist $H_{1},H_{2},H_{3}$ such that:

\emph{(i)} $I(H_{1};-1)=(-1)\cdot I(G;-1)$ and $\nu(H_{1})=\nu(G)$;

\emph{(ii)} $I(H_{2};-1)=I(G;-1)$ and $\nu(H_{2})=\nu(G)+1$;

\emph{(iii)} $\left\vert I(H_{3};-1)\right\vert =k\cdot I(G;-1)$ and
$\nu(H_{3})=\nu(G)+k-1$.
\end{corollary}

\begin{proof}
\emph{(i)} Let $H_{1}=H[v,K_{2},G]$ be the graph defined as in Lemma
\ref{lem5} (i.e., by joining with an edge the vertex $v$ to one of the
endpoints of $K_{2}$, say $u_{1}$, and to some vertex $u_{2}$ of $G$; see
Figure \ref{fig12}). \begin{figure}[h]
\setlength{\unitlength}{0.8cm} \begin{picture}(5,1.2)\thicklines
\multiput(4,0)(1,0){2}{\circle*{0.29}}
\put(4,0){\line(1,0){1}}
\put(4,0.4){\makebox(0,0){$a$}}
\put(5,0.4){\makebox(0,0){$u_{1}$}}
\put(3,0.5){\makebox(0,0){$K_{2}$}}
\multiput(8,0)(1,0){3}{\circle*{0.29}}
\put(12,0.5){\circle*{0.29}}
\put(8,0){\line(1,0){2}}
\put(10,0){\line(4,1){2}}
\put(11.5,0){\line(1,0){4}}
\put(11.5,0){\line(0,1){1}}
\put(11.5,1){\line(1,0){4}}
\put(15.5,0){\line(0,1){1}}
\put(8,0.4){\makebox(0,0){$a$}}
\put(9,0.4){\makebox(0,0){$u_{1}$}}
\put(10,0.4){\makebox(0,0){$v$}}
\put(12.5,0.5){\makebox(0,0){$u_{2}$}}
\put(14,0.5){\makebox(0,0){$G$}}
\put(7,0.5){\makebox(0,0){$H_{1}$}}
\end{picture}
\caption{The graphs $H_{1}$ and $G$ have the same cyclomatic number.}%
\label{fig12}%
\end{figure}
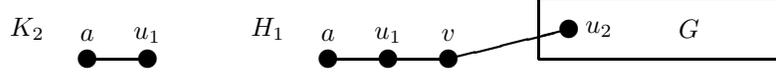

According to Lemma \ref{lem5} we get $I(H_{1};-1)=(-1)\cdot I(G;-1)$, since
$I(K_{2};x)=1+2x$ and $I(K_{2}-u_{1};x)=1+x$.

\emph{(ii)} Let $H_{2}=H[v,W,G]$ be the graph depicted in Figure \ref{fig232}.
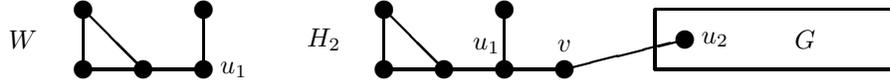
\begin{figure}[h]
\setlength{\unitlength}{0.8cm} \begin{picture}(5,1.2)\thicklines
\multiput(3,0)(1,0){3}{\circle*{0.29}}
\multiput(3,1)(2,0){2}{\circle*{0.29}}
\put(3,0){\line(1,0){2}}
\put(3,0){\line(0,1){1}}
\put(3,1){\line(1,-1){1}}
\put(5,0){\line(0,1){1}}
\put(5.5,0){\makebox(0,0){$u_{1}$}}
\put(2,0.5){\makebox(0,0){$W$}}
\multiput(8,0)(1,0){4}{\circle*{0.29}}
\multiput(8,1)(2,0){2}{\circle*{0.29}}
\put(13,0.5){\circle*{0.29}}
\put(8,0){\line(1,0){3}}
\put(8,0){\line(0,1){1}}
\put(8,1){\line(1,-1){1}}
\put(10,0){\line(0,1){1}}
\put(11,0){\line(4,1){2}}
\put(12.5,0){\line(1,0){4}}
\put(12.5,0){\line(0,1){1}}
\put(12.5,1){\line(1,0){4}}
\put(16.5,0){\line(0,1){1}}
\put(9.7,0.4){\makebox(0,0){$u_{1}$}}
\put(11,0.4){\makebox(0,0){$v$}}
\put(13.5,0.5){\makebox(0,0){$u_{2}$}}
\put(15,0.5){\makebox(0,0){$G$}}
\put(7,0.5){\makebox(0,0){$H_{2}$}}
\end{picture}
\caption{The graphs $H_{2}$ and $G$ satisfy $\nu(H_{2})=\nu(G)+1$.}%
\label{fig232}%
\end{figure}

Since $I(W;x)=1+5x+5x^{2}$ and $I(W-u_{1};x)=(1+x)(1+3x)$, Lemma \ref{lem5}
implies that $I(H_{2};-1)=I(G;-1)$.

\emph{(iii)} Let $L_{0}=K_{1}$ and $L_{s},s\geq1$, be the graph from Figure
\ref{fig212}. Notice that
\[
I(L_{0};x)=1,\quad I(L_{1};x)=1+3x,
\]
while, by Proposition \ref{prop1}\emph{(ii)}, $I(L_{s};x)$ satisfies
\[
I(L_{s};x)=(1+3x)\cdot I(L_{s-1};x)-x^{2}\cdot I(L_{s-2};x),s\geq2.
\]
Hence, it follows that $I(L_{s};-1)=(s+1)\cdot(-1)^{s},s\geq0$.
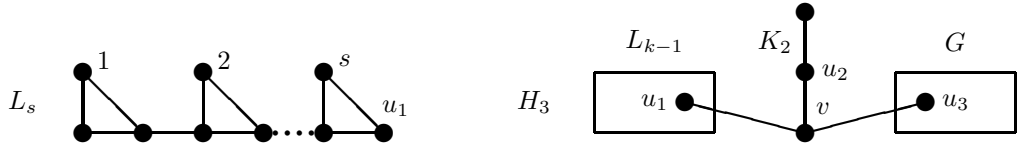
\begin{figure}[h]
\setlength{\unitlength}{0.8cm} \begin{picture}(5,2.2)\thicklines
\multiput(2,0)(1,0){6}{\circle*{0.29}}
\put(7.2,0.4){\makebox(0,0){$u_{1}$}}
\multiput(2,1)(2,0){3}{\circle*{0.29}}
\multiput(5,0)(0.2,0){5}{\circle*{0.1}}
\put(2,0){\line(1,0){3}}
\put(2,0){\line(0,1){1}}
\put(2,1){\line(1,-1){1}}
\put(4,1){\line(1,-1){1}}
\put(4,0){\line(0,1){1}}
\put(6,0){\line(1,0){1}}
\put(6,0){\line(0,1){1}}
\put(6,1){\line(1,-1){1}}
\put(2.35,1.2){\makebox(0,0){$1$}}
\put(4.35,1.2){\makebox(0,0){$2$}}
\put(6.35,1.2){\makebox(0,0){$s$}}
\put(1,0.5){\makebox(0,0){$L_{s} $}}
\put(12,0.5){\circle*{0.29}}
\multiput(14,1)(0,1){2}{\circle*{0.29}}
\put(16,0.5){\circle*{0.29}}
\put(14,0){\circle*{0.29}}
\put(14,0){\line(0,1){2}}
\put(14,0){\line(4,1){2}}
\put(14,0){\line(-4,1){2}}
\put(10.5,0){\line(1,0){2}}
\put(10.5,0){\line(0,1){1}}
\put(10.5,1){\line(1,0){2}}
\put(12.5,0){\line(0,1){1}}
\put(15.5,0){\line(1,0){2}}
\put(15.5,0){\line(0,1){1}}
\put(15.5,1){\line(1,0){2}}
\put(17.5,0){\line(0,1){1}}
\put(14.3,0.4){\makebox(0,0){$v$}}
\put(11.5,0.5){\makebox(0,0){$u_{1}$}}
\put(14.5,1){\makebox(0,0){$u_{2}$}}
\put(13.5,1.5){\makebox(0,0){$K_{2}$}}
\put(11.5,1.5){\makebox(0,0){$L_{k-1}$}}
\put(16.5,0.5){\makebox(0,0){$u_{3}$}}
\put(16.5,1.5){\makebox(0,0){$G$}}
\put(9.5,0.5){\makebox(0,0){$H_{3}$}}
\end{picture}
\caption{The graph $H_{3}$ has $\nu(H_{3})=\nu(G)+\nu(L_{k-1})=\nu(G)+k$.}%
\label{fig212}%
\end{figure}

Let $H_{3}=H[v,L_{k-1},K_{2},G]$ be the graph from Figure \ref{fig212}. Lemma
\ref{lem5} implies that
\[
I(H_{3};-1)=(-1)\cdot I(L_{k-1};-1)\cdot I(G;-1)=(-1)^{k}\cdot k\cdot
I(G;-1),
\]
which clearly ensures that $\left\vert I(H_{3};-1)\right\vert =k\cdot I(G;-1)
$.
\end{proof}

\begin{lemma}
\label{lem4}For every positive integer $\nu$, there exist connected graphs
$G_{1},G_{2}$, with $\nu(G_{1})=\nu(G_{2})=\nu$, such that $I(G_{1}%
;-1)=\pm2^{q}$ and $I(G_{2};-1)=\pm(2^{q}-1)$, for every $q\in\{0,1,2,...,\nu
\}$.
\end{lemma}

\begin{proof}
Let us remark that, by Corollary \ref{cor33}\emph{(i)}, to show that there is
a graph $G$ with $I(G;-1)=p$, for some integer $p$, it is enough to find a
graph $G$ satisfying $\left\vert I(G;-1)\right\vert =\left\vert p\right\vert $.

Let $G_{1}$ be the graph obtained from a $P_{3}=(\{u,v,w\},\{uv,vw\})$ and
joining $v$ to the endpoints of $\nu-q$ graphs isomorphic to $K_{2}$, while
$w$ is joined by an edge to one vertex of each of $q$ graphs isomorphic to
$K_{3}$ (see Figure \ref{fig911} for two examples). It is easy to see that
$\nu(G_{1})=\nu$. \begin{figure}[h]
\setlength{\unitlength}{0.8cm} \begin{picture}(5,3)\thicklines
\multiput(3,0)(0,1){3}{\circle*{0.29}}
\multiput(4,0)(0,1){4}{\circle*{0.29}}
\multiput(5,1)(0,1){3}{\circle*{0.29}}
\multiput(6,0)(0,1){4}{\circle*{0.29}}
\multiput(7,1)(0,2){2}{\circle*{0.29}}
\put(3,0){\line(1,0){1}}
\put(3,0){\line(1,1){1}}
\put(3,1){\line(1,0){4}}
\put(3,2){\line(1,0){1}}
\put(3,2){\line(1,-1){1}}
\put(4,0){\line(0,1){2}}
\put(4,3){\line(1,0){1}}
\put(4,3){\line(1,-1){1}}
\put(5,1){\line(0,1){2}}
\put(5,1){\line(1,1){2}}
\put(6,0){\line(1,1){1}}
\put(6,0){\line(0,1){1}}
\put(6,2){\line(0,1){1}}
\put(6,3){\line(1,0){1}}
\put(3,0.6){\makebox(0,0){$u$}}
\put(4.2,0.6){\makebox(0,0){$v$}}
\put(5.1,0.6){\makebox(0,0){$w$}}
\put(2,1.5){\makebox(0,0){$H_{1}$}}
\multiput(10,1)(1,0){5}{\circle*{0.29}}
\multiput(10,2)(1,0){6}{\circle*{0.29}}
\multiput(11,3)(1,0){5}{\circle*{0.29}}
\put(14,0){\circle*{0.29}}
\put(10,1){\line(1,0){4}}
\put(10,2){\line(1,-1){1}}
\put(10,2){\line(1,0){1}}
\put(11,1){\line(0,1){1}}
\put(11,3){\line(1,0){1}}
\put(11,3){\line(1,-1){1}}
\put(12,1){\line(1,1){2}}
\put(12,1){\line(0,1){2}}
\put(12,1){\line(2,1){2}}
\put(13,1){\line(1,-1){1}}
\put(13,2){\line(0,1){1}}
\put(13,3){\line(1,0){1}}
\put(14,2){\line(1,1){1}}
\put(14,2){\line(1,0){1}}
\put(14,0){\line(0,1){1}}
\put(15,2){\line(0,1){1}}
\put(10,0.6){\makebox(0,0){$u$}}
\put(11,0.6){\makebox(0,0){$v$}}
\put(12,0.6){\makebox(0,0){$w$}}
\put(9,1.5){\makebox(0,0){$H_{2}$}}
\end{picture}
\caption{$\nu(H_{1})=\nu(H_{2})=5$, while $I(H_{1};-1)=(-1)^{4}\cdot2^{3}$ and
$I(H_{2};-1)=(-1)^{5}\cdot2^{4}$.}%
\label{fig911}%
\end{figure}
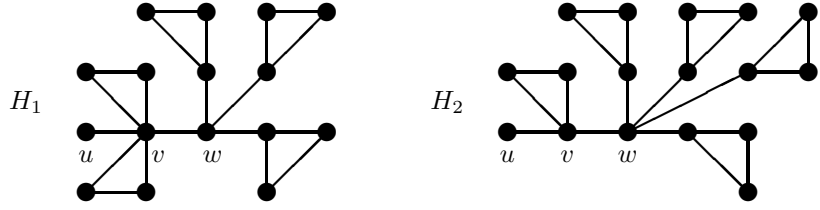

Notice that $u$ is a pendant vertex of $G_{1}$ and $v\in N(u)$. Hence, using
Lemma \ref{lem2}, we infer that
\begin{align*}
I(G_{1};-1)  &  =(-1)\cdot I(G_{1}-N[v];-1)\\
&  =(-1)\cdot I(qK_{3};-1)=(-1)^{q+1}\cdot2^{q}.
\end{align*}

Let $W_{q},q\geq2$, be the graph obtained from the star $K_{1,q}$ by adding
$qK_{2}$ , such that each pendant vertex of $K_{1,q}$ is joined to the
endpoints of one $K_{2}$ (see Figure \ref{fig88} for $W_{5}$).
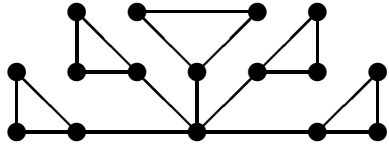
\begin{figure}[h]
\setlength{\unitlength}{0.8cm} \begin{picture}(5,3)\thicklines
\multiput(6,0.5)(3,0){3}{\circle*{0.29}}
\multiput(6,1.5)(1,0){7}{\circle*{0.29}}
\multiput(7,2.5)(1,0){2}{\circle*{0.29}}
\multiput(10,2.5)(1,0){2}{\circle*{0.29}}
\put(7,0.5){\circle*{0.29}}
\put(11,0.5){\circle*{0.29}}
\multiput(6,0.5)(1,0){6}{\line(1,0){1}}
\multiput(6,0.5)(3,0){3}{\line(0,1){1}}
\put(7,1.5){\line(0,1){1}}
\put(11,1.5){\line(0,1){1}}
\put(8,1.5){\line(-1,1){1}}
\put(10,1.5){\line(1,1){1}}
\put(7,1.5){\line(1,0){1}}
\put(9,0.5){\line(-1,1){1}}
\put(9,0.5){\line(1,1){1}}
\put(10,1.5){\line(1,0){1}}
\put(7,0.5){\line(-1,1){1}}
\put(11,0.5){\line(1,1){1}}
\put(9,1.5){\line(-1,1){1}}
\put(9,1.5){\line(1,1){1}}
\put(8,2.5){\line(1,0){2}}
\end{picture}
\caption{The graph $W_{5}$ has $I(W_{5};x)=(1+3x)^{5}+x\cdot(1+2x)^{5}$ and
hence, $I(W_{5};-1)=(-1)^{5}(2^{5}-1)$.}%
\label{fig88}%
\end{figure}

Then, $I(W_{q};x)=(1+3x)^{q}+x\cdot(1+2x)^{q}$ and consequently, it follows
that
\[
I(W_{q};-1)=(-2)^{q}-(-1)^{q}=(-1)^{q}(2^{q}-1).
\]

Using Corollary \ref{cor33}\emph{(ii)} for $\nu-q$ times one can build a graph
$G_{2}$ such that $\nu(G_{2})=\nu$, while $I(G_{2};-1)=I(W_{q};-1)=(-1)^{q}%
(2^{q}-1)$.
\end{proof}

\begin{proposition}
For any positive integer $q=p_{1}^{\varepsilon_{1}}\cdot...p_{k}%
^{\varepsilon_{k}}$, where $p_{1},p_{2},...,p_{k}$ are prime numbers, there is
a graph $G$ such that
\[
\nu(G)=\varepsilon_{1}(p_{1}-1)+...+\varepsilon_{k}(p_{k}-1)\text{ and
}I(G;-1)=q.
\]

\end{proposition}

\begin{proof}
Recall that the graph $L_{s}$, depicted in Figure \ref{fig212}, has
$I(L_{s};-1)=(-1)^{s}\cdot(s+1)$.

The graph $G$ is obtained starting with a $P_{3}=(\{u,v,w\},\{uv,vw\})$ and
joining $w$ to one vertex from each of $\varepsilon_{j}$ graphs isomorphic to
$L_{p_{j}-1}$ (see Figure \ref{fig11} for an example). Firstly, it is easy to
see that $\nu(G)=\varepsilon_{1}(p_{1}-1)+...+\varepsilon_{k}(p_{k}-1)$.
\begin{figure}[h]
\setlength{\unitlength}{0.8cm} \begin{picture}(5,2.3)\thicklines
\multiput(5,0)(1,0){9}{\circle*{0.29}}
\multiput(5,1)(1,0){8}{\circle*{0.29}}
\multiput(5,2)(1,0){7}{\circle*{0.29}}
\put(5,0){\line(1,0){8}}
\put(5,1){\line(0,1){1}}
\put(5,2){\line(1,-1){2}}
\put(5,1){\line(1,0){1}}
\put(6,2){\line(1,0){1}}
\put(6,2){\line(1,-1){1}}
\put(7,0){\line(1,1){2}}
\put(7,0){\line(0,1){2}}
\put(8,0){\line(1,1){1}}
\put(8,1){\line(0,1){1}}
\put(8,2){\line(1,0){3}}
\put(9,0){\line(0,1){1}}
\put(9,1){\line(1,-1){1}}
\put(10,0){\line(0,1){1}}
\put(10,2){\line(1,-1){1}}
\put(11,1){\line(0,1){1}}
\put(12,1){\line(1,-1){1}}
\put(12,0){\line(0,1){1}}
\put(7,-0.4){\makebox(0,0){$w$}}
\put(4,1){\makebox(0,0){$G$}}
\end{picture}
\caption{The vertex $w$ of $G$ is joined to two $L_{1}$, one $L_{2}$ and one
$L_{3}$.}%
\label{fig11}%
\end{figure}
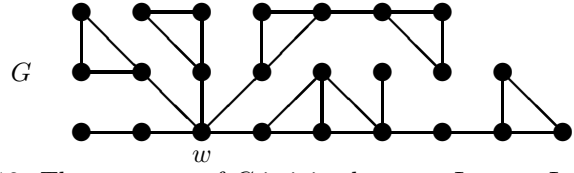

Then, applying Proposition \ref{prop1}\emph{(i)}, we get
\begin{align*}
I(G;x)  &  =I(G-w;x)+x\cdot I(G-N[w];x)\\
&  =(1+2x)\cdot\left(  I(L_{p_{1}-1};x)\right)  ^{\varepsilon_{1}}%
\cdot...\cdot\left(  I(L_{p_{k}-1};x)\right)  ^{\varepsilon_{k}}\\
&  +x\cdot(1+x)\cdot I(G-N[w]\cup\{u\};x),
\end{align*}
which implies that $\left\vert I(G;-1)\right\vert =p_{1}^{\varepsilon_{1}%
}\cdot...\cdot p_{k}^{\varepsilon_{k}}$.
\end{proof}

\section{Conclusions}

In this paper we\emph{\ }proved that $\left\vert I(G;-1)\right\vert \leq
2^{\nu(G)}$\ and we exhibit graphs satisfying $\left\vert I(G;-1)\right\vert
\in\{2^{\nu(G)},2^{\nu(G)}-1\}$.

\begin{conjecture}
For every positive integer $\nu$ and each integer $q$ such that $\left\vert
q\right\vert \leq2^{\nu}$, there is a connected graph $G$ with $\nu(G)=\nu$
and $I(G;-1)=q$.
\end{conjecture}

Taking into account Corollary \ref{cor33}\emph{(i)}, it is enough to find a
connected graph $G$ such that $\left\vert I(G;-1)\right\vert =\left\vert
q\right\vert $.

\begin{remark}
We have the complete solution for $\nu\in\{0,1,2,3\}$.
\end{remark}

\begin{itemize}
\item If $\nu(G)=0$, then $G$ must be a forest and for every $q\in
\{-2^{0},0,2^{0}\}=\{-1,0,1\}$, we saw above that there is a tree $T$ such
that $I(T;-1)=q$.

\item If $\nu(G)=1$, then $G$ must be a unicycle graph and for every
$q\in\{-2,-1,0,1,2\}$, there is a connected graph $G$ such that $I(G;-1)=q$.

\item Let $\nu=2$ and $q\in\{-2^{2},-3,...,3,2^{2}\}$. According to Corollary
\ref{cor33}\emph{(ii)},\emph{(iii)}, and Lemma \ref{lem4}, it follows that
there is a connected graph $G$ such that $I(G;-1)=q$.

\item Let $\nu=3$ and $q\in\{-2^{3},-7,...,7,2^{3}\}$. By Corollary
\ref{cor33}\emph{(ii)} and \emph{(iii)}, and Lemma \ref{lem4}, it follows that
there is a connected graph $G$ such that
\[
I(G;-1)=q\in\{-2^{3},-7,-6,-4,-3,-2,-1,0,1,2,3,4,6,7,2^{3}\}.
\]
The graph $G$ from Figure \ref{fig22} has
\begin{align*}
I(G;x)  &  =I(G-v;x)+x\cdot I(G-N[v];x)\\
&  =(1+2x)\cdot(1+3x)^{2}\cdot I(P_{3};x)+x\cdot(1+2x)^{3},
\end{align*}
which implies $I(G;-1)=5$. \begin{figure}[h]
\setlength{\unitlength}{0.8cm} \begin{picture}(5,2.3)\thicklines
\multiput(3,0)(1,0){5}{\circle*{0.29}}
\multiput(3,1)(1,0){5}{\circle*{0.29}}
\multiput(6,2)(1,0){2}{\circle*{0.29}}
\put(3,0){\line(1,0){4}}
\put(3,0){\line(0,1){1}}
\put(4,1){\line(1,-1){1}}
\put(4,1){\line(1,0){1}}
\put(5,0){\line(0,1){1}}
\put(5,0){\line(1,1){2}}
\put(6,2){\line(1,0){1}}
\put(6,1){\line(0,1){1}}
\put(6,0){\line(1,1){1}}
\put(7,0){\line(0,1){1}}
\put(5.5,0.2){\makebox(0,0){$v$}}
\put(2,1){\makebox(0,0){$G$}}
\multiput(10,0)(1,0){7}{\circle*{0.29}}
\multiput(11,1)(1,0){5}{\circle*{0.29}}
\multiput(11,2)(1,0){3}{\circle*{0.29}}
\put(10,0){\line(1,0){6}}
\put(11,1){\line(0,1){1}}
\put(11,1){\line(1,0){1}}
\put(11,2){\line(1,-1){1}}
\put(12,0){\line(0,1){1}}
\put(12,0){\line(1,1){1}}
\put(12,2){\line(1,0){1}}
\put(12,2){\line(1,-1){1}}
\put(13,0){\line(1,1){1}}
\put(13,1){\line(0,1){1}}
\put(14,0){\line(0,1){1}}
\put(15,1){\line(1,-1){1}}
\put(15,0){\line(0,1){1}}
\put(11.7,0.3){\makebox(0,0){$v$}}
\put(9,1){\makebox(0,0){$H$}}
\end{picture}
\caption{The graphs $G$ and $H$ satisfy $\nu(G)=3$ and $\nu(H)=4$.}%
\label{fig22}%
\end{figure}
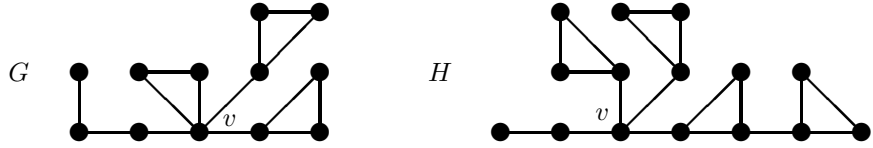

\item Let $\nu=4$ and $q\in A=\{-2^{4},-15,...,15,2^{4}\}$. By Corollary
\ref{cor33}\emph{(ii)} and \emph{(iii)}, and Lemma \ref{lem4}, it follows that
there is a connected graph $G$ such that
\[
I(G;-1)=q\in A-\{-13,-11,11,13\}.
\]
The graph $H$ from Figure \ref{fig22} has
\begin{align*}
I(H;x)  &  =I(H-v;x)+x\cdot I(H-N[v];x)\\
&  =(1+3x)^{2}\cdot I(P_{3};x)\cdot I(\bigtriangleup_{2};x)+x\cdot
(1+2x)^{3}\cdot\left(  1+5x+5x^{2}\right)  ,
\end{align*}
which implies $I(H;-1)=11$. The case of $q=13$ has not been settled yet.
\end{itemize}

\end{document}